\documentclass[12pt]{article} 
\usepackage{color} 
\usepackage{amssymb} 
\usepackage{amsmath} 
\usepackage{graphicx}
\usepackage{multirow} 
\usepackage{longtable} 
\usepackage{rotating} 
\usepackage{longtable} 
\usepackage{lscape} 
\usepackage[font=scriptsize]{caption} 
\date{} 
\title{On the Parameter Estimation of the Generalized Exponential Distribution Under Progressive \\ Type-I Interval Censoring Scheme\\
\vspace{1cm}
\small{{Mahdi Teimouri}}\\
Email: teimouri@aut.ac.ir\\
Department of Mathematics and Statistics, Faculty of Science and Engineering, \\
Gonbad Kavous University, Gonbad Kavous, Iran.} 
\begin{document}
\maketitle{} 
\noindent{} {\bf{Abstract:}} Chen and Lio (Computational Statistics and Data Analysis 54: 1581-1591, 2010) proposed five methods for 
estimating the parameters of generalized exponential distribution under progressive type-I interval censoring scheme. Unfortunately, among them, the proposed EM algorithm is incorrect. Here, we propose the correct EM algorithm and compare its performance with the maximum likelihood estimators and that proposed by Chen and Lio (2010) in a simulation study.
\\\\
\noindent{} {\bf{Keyword:}} EM algorithm; Generalized exponential distribution; Maximum likelihood estimator; Progressive type-I interval censoring scheme.
\section{Introduction}
\subsection{Generalized exponential (GE) distribution}
The random variable $X$ follows GE distribution if its probability density function (pdf) and distribution function are given by
\begin{align}\label{pdf}
f(x,\theta)=\alpha \bigl(1-e^{-\lambda x}\bigr)^{\alpha-1} e^{-\lambda x},
\end{align}
and
\begin{align}\label{cdf}
F(x,\theta)=\bigl(1-e^{-\lambda x}\bigr)^{\alpha},
\end{align}
where $\theta=(\alpha,\lambda)$ is parameter vector ($\alpha$ is the shape parameter and $\lambda$ is the rate parameter). The family of GE distributions was introduced by Mudholkar and Srivastava (1993). For a comprehensive account of the theory and applications of GE distribution, we refer the readers to Gupta and Kundu (2007).
\subsection{Progressively type-I interval censoring scheme}
Suppose $n$ subjects are placed on a life testing simultaneously at time $t_0=0$ and under inspection at $m$ pre-determined times $t_1<t_2<\dots< t_m$ in which $t_m$ is the time to terminate the life testing. At the $i$-th inspection time, $t_i$, the number, $X_i$,
of failures within $(t_i, t_{i+1}]$ is recorded and $R_i$ alive items are randomly removed from the life testing, for $i=1,\dots,m$.
As pointed out by Chen and Lio (2010), since the number, $Y_i$, of surviving items is a random variable and the exact number of items withdrawn should not be greater than $Y_i$ at time schedule $t_i$, then $R_i$ could be determined by the pre-specified percentage of the remaining surviving units at $t_i$, or equivalently $R=\lfloor p_i Y_i\rfloor$; for $i=1,\dots,m$. Each progressively type-I interval censoring scheme is shown by $\{X_i,R_i,T_i\}_{i=1}^{m}$ where $n=\sum_{i=1}^{m}X_i+R_i$ is the sample size. If $Ri =0$; for $i=1,\dots,m-1$, then the progressively type-I interval censoring scheme is equivalent to a type-I interval censoring scheme with sample $X_1, X2, \dots,X_m, X_{m-1}=R_m$. Suppose a $\{X_i,R_i,T_i\}_{i=1}^{m}$ life testing scheme where $n$ items each follows independently the cdf $F(.,\theta)$ is under the test. The likelihood function is (see \cite{Aggarwala2001}) is 
\begin{align}\label{lik}
L(\theta)\propto \prod_{i=1}^{m}\bigl[F(t_{i},\theta)-F(t_{i-1},\theta)\bigr]^{X_{i}}\bigl[1-F(t_{i},\theta)\bigr]^{R_{i}}.
\end{align}
As the most common used tool, the maximum likelihood (ML) approach is employed to estimate the $\theta$. But, equation (\ref{lik}) must be maximized through iterative algorithm such as Newton-Raphson to obtain the ML estimators and there is no guarantee that the Newton-Raphson method converges. Another technique is the expectation-maximization (EM) algorithm that always converges, see \cite{Little}. However, if practitioner is interested in the ML estimators, the first few steps of the EM algorithm can be used to get a good starting value for the Newton-Raphson algorithm, see \cite{Nadarajah2014}.
\subsection{EM algorithm}\label{em}
The EM algorithm, introduced by \cite{Dempster1977}, is known as the popular method for computing the ML estimators when we encounter the incomplete data problem. In other word, the use of the EM algorithm involves cases that we are dealing with the latent variables, provided that the statistical model is formulated as a missing or latent variable problem. In what follows, we give a brief description of the EM algorithm. Let ${\boldsymbol{\xi}}$, ${\boldsymbol{Z}}$, and ${\boldsymbol{\omega}}$ denote the complete, unobservable variable, and observed data, respectively (complete data consists of observed values and unobservable variables, i.e., ${\boldsymbol{\xi}}=(\boldsymbol{Z},\boldsymbol{\omega})$). The EM algorithm works by maximizing the conditional expectation $Q\left(\theta |\theta^{(t)}\right)=E\left(l_c(\theta;\boldsymbol{\xi})|\boldsymbol{\omega}, \theta^{(t)}\right)$ of complete data log-likelihood function given observed data and a current estimate $\theta^{(t)}$ of the parameter vector $\theta$ where $l_c(\theta;\boldsymbol{x})$ denotes the complete data log-likelihood function. Each iteration of the EM algorithm consists of two steps:
\begin {enumerate} 
\item Expectation (E)-step: Computing $Q\left(\theta|\theta^{(t)}\right)$ at the $t$-th iteration. 
\item Maximization (M)-step: Maximizing $Q\left(\theta |\theta^{(t)}\right)$ with respect to $\theta$ to get $\theta^{(t+1)}$. 
\end {enumerate} 
The E-step and M-step are repeated until convergence occurs, see \cite{Dempster1977} and \cite{McLachlan1997}.
\section{EM algorithm for GE family under progressive type-I interval censoring scheme}
Suppose $n$ failure times follow the GE distribution with pdf and cdf given by expressions (\ref{pdf}) and (\ref{cdf}), respectively. For convenience, let us to use the notations given by Chen and Lio (2010). So, let $t_{i,j}$; for $j=1,\dots, X_i$, denote the independent and identically distributed (iid) failure times in the subinterval $(t_{i-1}, t_i]$; for $i=1,\dots, m$ and $t^{*}_{I,j}$; for $j=1,\dots,R_{i}$, indicate on iid failure times of the randomly removed items alive at the end of the subinterval $(t_{i-1}, t_i]$; for $i=1,\dots, m$. Then, the complete data log-likelihood, $l_{c}(\theta)$, is (see \cite{Chen2010})
\begin{align}\label{loglik0}
l_{c}(\theta)\propto \sum_{i=1}^{m}\sum_{j=1}^{X_{i}}\log f(T_{i,j},\theta)+\sum_{i=1}^{m}\sum_{j=1}^{R_{i}}\log f(T^{*}_{i,j},\theta).
\end{align}
In expression (\ref{loglik0}), we show unobservable (or missing) variables by capital letters $T_{i,j}$ (for $i=1,\dots, m$, $j=1,\dots,x_{i}$) and $T^{*}_{i,j}$ (for $i=1,\dots, m$; $j=1,\dots,r_{i}$) in which $\sum_{i=1}^{m}x_i+r_i=n$. The progressive type-I censoring scheme is an incomplete data problem. The observed values are $x_i$s and $r_i$s; for $i=1,\dots, m$ and unobservable variables are $T_{i,j}$ (iid failure times during subinterval $(t_{i-1}, t_i]$) and $T^{*}_{i,j}$ (iid withdrawn survival times during subinterval $(t_{i-1}, t_i]$). Therefore, under EM algorithm framework mentioned in subsection \ref{em}, the vector of observed data, $\boldsymbol{\omega}$, is $\boldsymbol{\omega}=(x_1,\dots,x_m,r_1,\dots,r_m)$ and the vector of unobservable variables is, $\boldsymbol{Z}=(T_{i,1},\dots,T_{i,x_i},T^{*}_{i,1},\dots,T^{*}_{i,r_i})$; for $i=1,\dots, m$. Assuming that we are at $t$-th iteration, in order to implement the EM algorithm, we follow two steps given by the following.
\begin{itemize}
\item {\bf{E-step}}: we need to compute the conditional expectation $Q\bigl(\theta \big|\theta^{(t)}\bigr)=E\bigl(l_c(\theta;\boldsymbol{\xi})\big|\boldsymbol{\omega}, \theta^{(t)}\bigr)$ of the complete data log-likelihood function. It follows, form (\ref{loglik0}), that
\begin{align}\label{loglik1}
Q\bigl(\theta \big|\theta^{(t)}\bigr)=&\text{C}+ \sum_{i=1}^{m}E \Bigl(\sum_{j=1}^{X_{i}}\log f(T_{i,j},\theta)\Big |\boldsymbol{\omega}, \theta=\theta^{(t)}\Bigr)+\sum_{i=1}^{m}E \Bigl(\sum_{j=1}^{R_{i}} \log f(T^{*}_{i,j},\theta) \Big|\boldsymbol{\omega}, \theta=\theta^{(t)}\Bigr)\nonumber\\
=&\text{C}+ \sum_{i=1}^{m}E \Bigl(\sum_{j=1}^{X_{i}}\log f(T_{i,j},\theta)\Big |X_i=x_i,T_{i,j} \in (t_{i-1}, t_i], \theta=\theta^{(t)}\Bigr)\nonumber\\
&+\sum_{i=1}^{m}E \Bigl(\sum_{j=1}^{R_{i}} \log f(T^{*}_{i,j},\theta) \Big|R_i=r_i,T^{*}_{i,j} \in [t_i,\infty), \theta=\theta^{(t)}\Bigr)\nonumber\\
=&\text{C}+ \sum_{i=1}^{m}\sum_{j=1}^{x_{i}}E \bigl(\log f(T_{i,j},\theta)\big |T_{i,j} \in (t_{i-1}, t_i], \theta=\theta^{(t)}\bigr)\nonumber\\
&+\sum_{i=1}^{m}\sum_{j=1}^{r_{i}}E \bigl(\log f(T^{*}_{i,j},\theta) \big|T^{*}_{i,j} \in [t_i,\infty), \theta=\theta^{(t)}\bigr)\nonumber\\
=&\bigl(\log (\alpha) +\log (\lambda)\bigr)\Bigl(\sum_{i=1}^{m}(x_i+r_i)\Bigr)
-\lambda \sum_{i=1}^{m}\sum_{j=1}^{x_i}E \bigl(T_{i,j} \big |T_{i,j} \in (t_{i-1}, t_i],, \theta=\theta^{(t)}\bigr)\nonumber\\
&-\lambda \sum_{i=1}^{m}\sum_{j=1}^{r_i}E \bigl(T^{*}_{i,j} \big |T^{*}_{i,j} \in [t_i,\infty), \theta=\theta^{(t)}\bigr)\nonumber\\
&+(\alpha-1) \sum_{i=1}^{m}\sum_{j=1}^{x_i}E \bigl(\log \bigl(1-e^{-\lambda T_{i,j}} \bigr)\big |T_{i,j} \in (t_{i-1}, t_i], \theta=\theta^{(t)}\bigr)\nonumber\\
&+(\alpha-1) \sum_{i=1}^{m}\sum_{j=1}^{r_i}E \bigl(\log \bigl(1-e^{-\lambda T^{*}_{i,j}} \bigr)\big |T^{*}_{i,j} \in [t_i,\infty), \theta=\theta^{(t)}\bigr),
\end{align}
where $\text{C}$ is a constant independent of $\theta$ and $\theta^{(t)}=(\alpha^{(t)},\lambda^{(t)})$. We note that the lifetimes of the $r_i$ unobserved items during subinterval $(t_{i-1}, t_i]$ are conditionally independent, identically distributed, and follow the truncated GE distribution on interval $[t_{i}, \infty)$. Also, lifetimes of the $x_i$ unobservable subjects during subinterval $(t_{i-1}, t_i]$ are conditionally independent, identically distributed, and follow the double-truncated GE distribution on subinterval $(t_{i-1}, t_i]$; for $i=1,\dots, m$. Therefore, considering the right-hand side of (\ref{loglik1}), the required conditional expectations are:
\begin{align}
E_{1i}&=E \Bigl(T_{i,j} \Big | T_{i,j} \in (t_{i-1}, t_i], \theta^{(t)}=(\alpha^{(t)},\lambda^{(t)})\Bigr)=\frac{\int_{t_{i-1}}^{t_{i}}uf\bigl(u,\theta^{(t)}\bigr)du}{F\bigl(t_{i},\theta^{(t)}\bigr)-F\bigl(t_{i-1},\theta^{(t)}\bigr)}\label{e1i},\\
E_{2i}&=E \Bigl(\log \bigl(1-e^{-\lambda T_{i,j}} \bigr) \Big | T_{i,j} \in (t_{i-1}, t_i], \theta^{(t)}=(\alpha^{(t)},\lambda^{(t)})\Bigr)\nonumber\\
&=\frac{\int_{t_{i-1}}^{t_{i}}\log \bigl(1-e^{-\lambda u} \bigr)f\bigl(u,\theta^{(t)}\bigr)du}{F\bigl(t_{i},\theta^{(t)}\bigr)-F\bigl(t_{i-1},\theta^{(t)}\bigr)},\label{e2i}\\
E_{3i}&=E \Bigl(T^{*}_{i,j} \Big | T^{*}_{i,j} \in [t_i,\infty), \theta^{(t)}=(\alpha^{(t)},\lambda^{(t)})\Bigr)=\frac{\int_{t_{i}}^{\infty}uf\bigl(u,\theta^{(t)}\bigr)du}{1-F\bigl(t_{i},\theta^{(t)}\bigr)},\label{e3i}\\
E_{4i}&=E \Bigl(\log \bigl(1-e^{-\lambda T^{*}_{i,j}} \bigr) \big | T^{*}_{i,j} \in [t_{i}, \infty), \theta^{(t)}=(\alpha^{(t)},\lambda^{(t)})\Bigr)\nonumber\\
&=\frac{\int_{t_{i}}^{\infty}\log \bigl(1-e^{-\lambda u} \bigr)f\bigl(u,\theta^{(t)}\bigr)du}{1-F\bigl(t_{i},\theta^{(t)}\bigr)},\label{e4i}
\end{align}
where $i=1,\dots, m$ and $t_0=0$.
\item {\bf{M-step}}: by substituting the computed conditional expectations $E_{1i}, E_{2i}, E_{3i}$, and $E_{4i}$ given in (\ref{e1i})-(\ref{e4i}) into the right-hand side of (\ref{loglik1}), we follow the EM algorithm by calculating the derivatives with respect to parameters as follows.
\begin{align}
\frac{\partial Q\bigl(\theta \big|\theta^{(t)}\bigr)}{\partial \alpha}=\frac{\sum_{i=1}^{m}(x_i+r_i)}{\alpha}+\sum_{i=1}^{m} x_iE_{2i}+\sum_{i=1}^{m} r_i E_{4i},\label{alpha0}\\ 
\frac{\partial Q\bigl(\theta \big|\theta^{(t)}\bigr)}{\partial \lambda}=\frac{\sum_{i=1}^{m}(x_i+r_i)}{\lambda}-\sum_{i=1}^{m} x_i E_{1i}-\sum_{i=1}^{m} r_iE_{3i},\label{lambda0}
\end{align}
where $\sum_{i=1}^{m}(x_i+r_i)=n$. Equating the right-hand side of (\ref{alpha0}) and (\ref{lambda0}) to zero it turns out that
\begin{align}
\alpha^{(t)}=-\frac{n}{\sum_{i=1}^{m} x_i E_{2i}+\sum_{i=1}^{m} r_i E_{4i}},\label{alphahat}
\end{align}
and
\begin{align}
\lambda^{(t)}=\frac{n}{\sum_{i=1}^{m} x_i E_{1i}+\sum_{i=1}^{m} r_i E_{3i}}.\label{lambdahat}
\end{align}
The M-step is complete.
\end{itemize}
We mention that the EM algorithm proposed by Chen and Lio (2010) is incorrect since they took expectation form the complete data log-likelihood function after differentiating it with respect to parameters which in not usual in the EM framework. Using the starting values as $\theta^{(0)}=(\alpha^{(0)},\lambda^{(0)})$ and repeating the E-step and M-step described as above the EM estimators are obtained. Compare the updated shape and rate parameters at $t$-th iteration given in (\ref{alphahat}) and (\ref{lambdahat}) with those given by Chen and Lio (2010). It is known that the updated shape parameters are the same but there is a significant difference between updated rate parameter given here and that given in Chen and Lio (2010). Although, difference between rate parameters is theoretically significant, however we perform a simulation study in the next section to observe the differences visually. 
\section{Simulation study}
\begin{figure}
\resizebox{\textwidth}{!}
{\begin{tabular}{ccc}
\includegraphics[width=50mm,height=50mm]{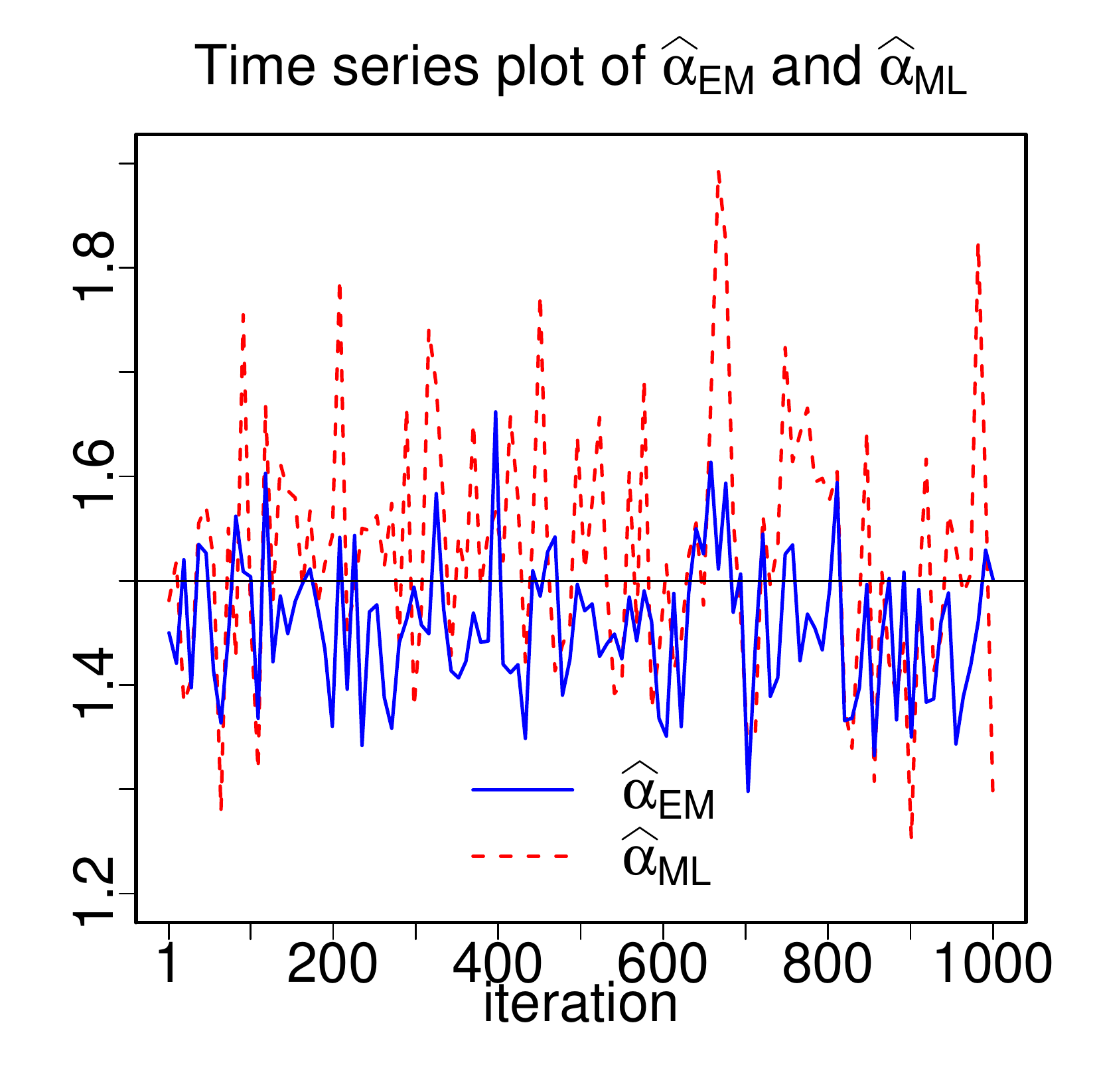}&
\includegraphics[width=50mm,height=50mm]{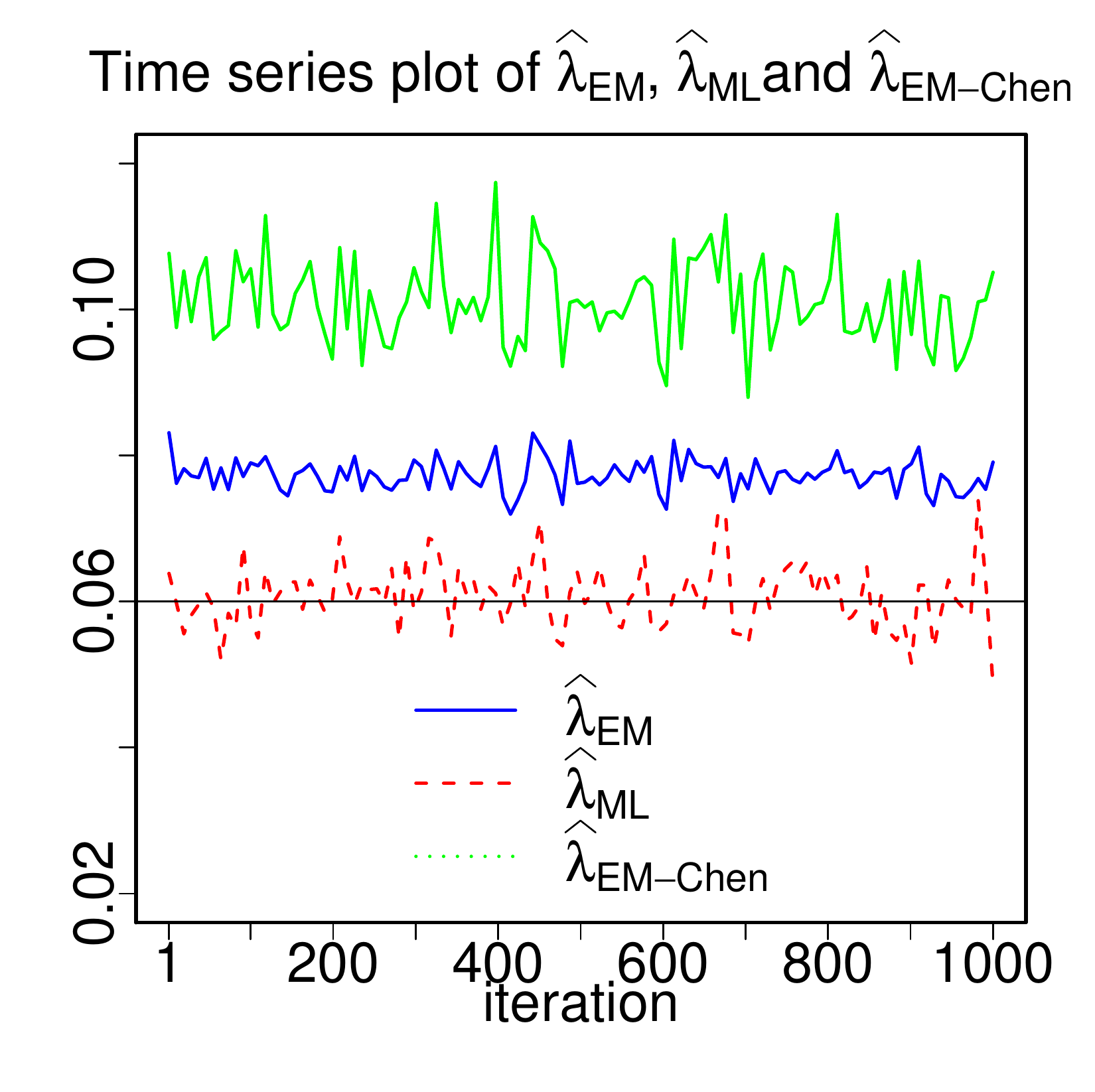}\\
\includegraphics[width=50mm,height=50mm]{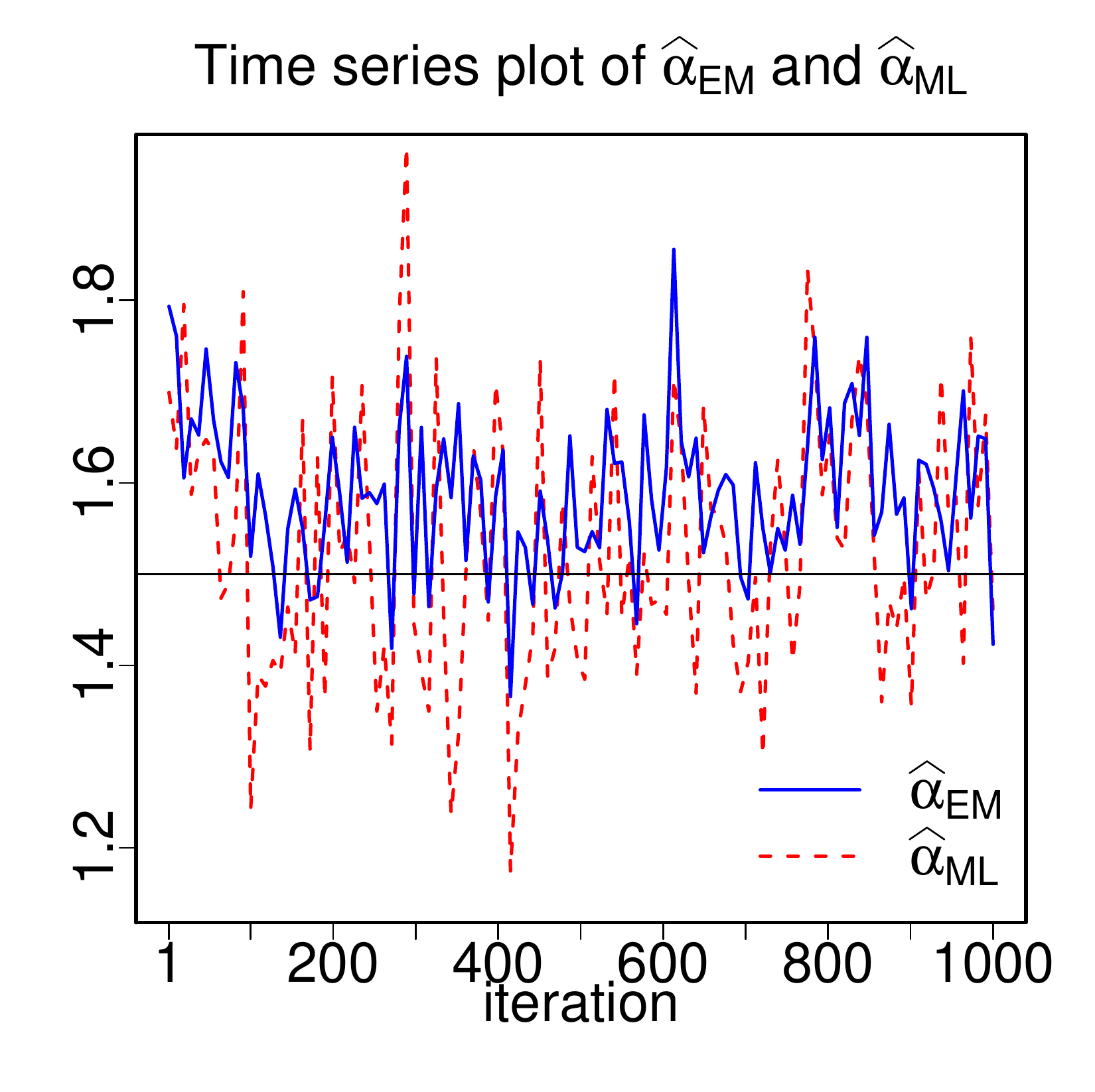}&
\includegraphics[width=50mm,height=50mm]{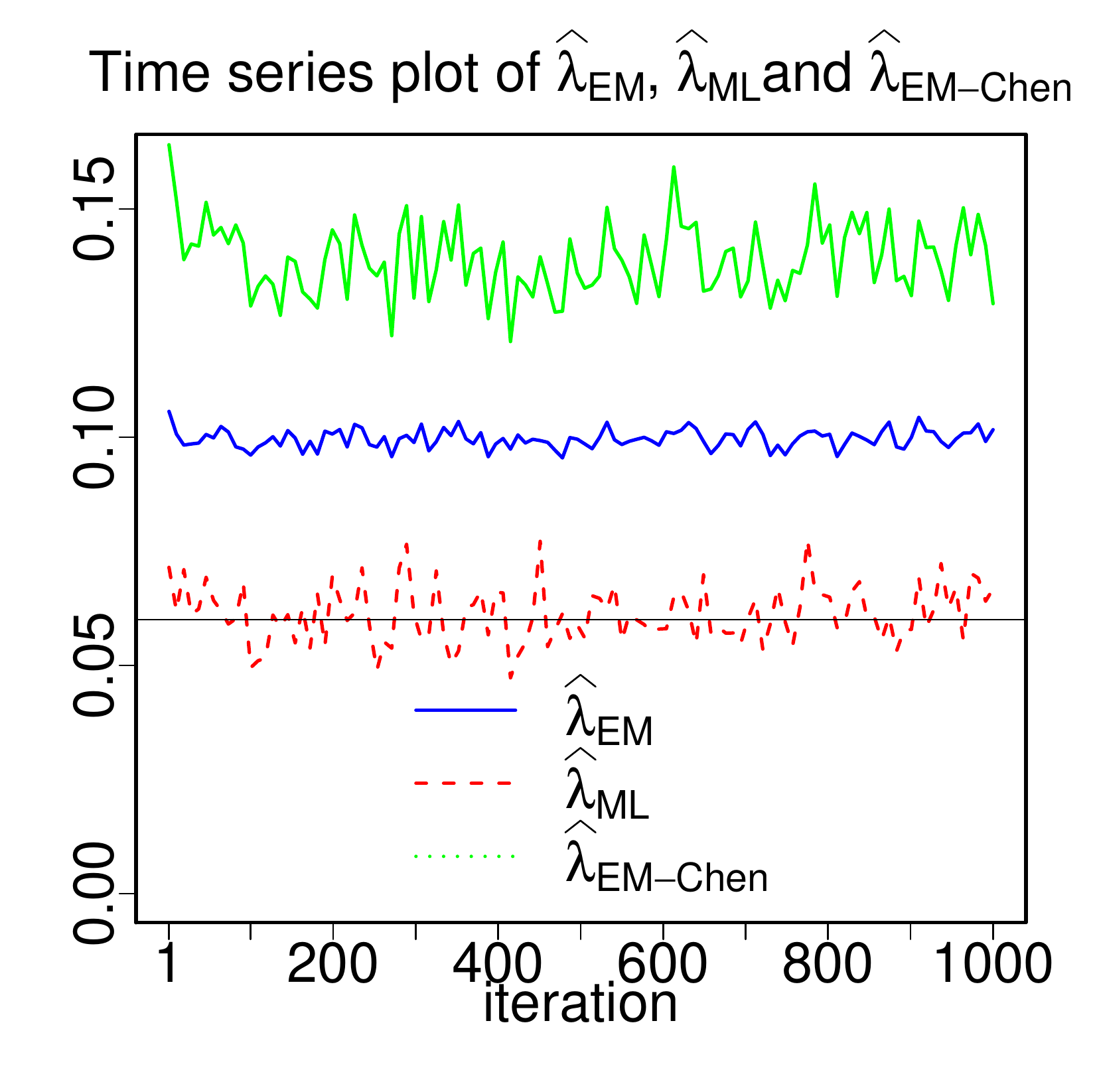}\\
\end{tabular}}
\caption{Time series plot of $\hat{\alpha}_{EM}$, $\hat{\alpha}_{ML}$, and $\hat{\alpha}_{EM-Chen}$ under settings $p_{(1)}$ (top row) and $p_{(2)}$ (bottom row).}
\label{fig1}
\end{figure}

\begin{figure}
\resizebox{\textwidth}{!}
{\begin{tabular}{ccc}
\includegraphics[width=50mm,height=50mm]{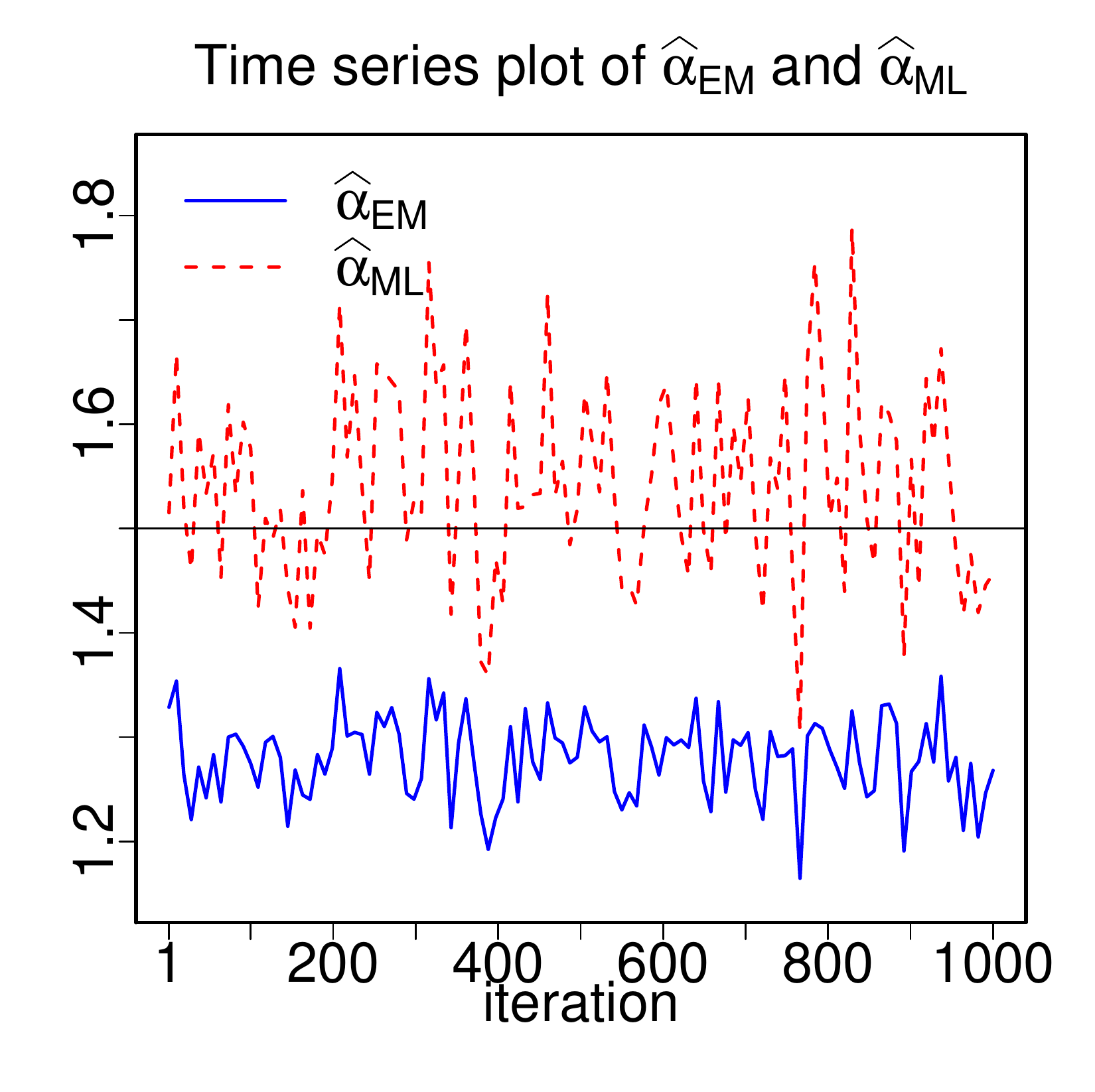}&
\includegraphics[width=50mm,height=50mm]{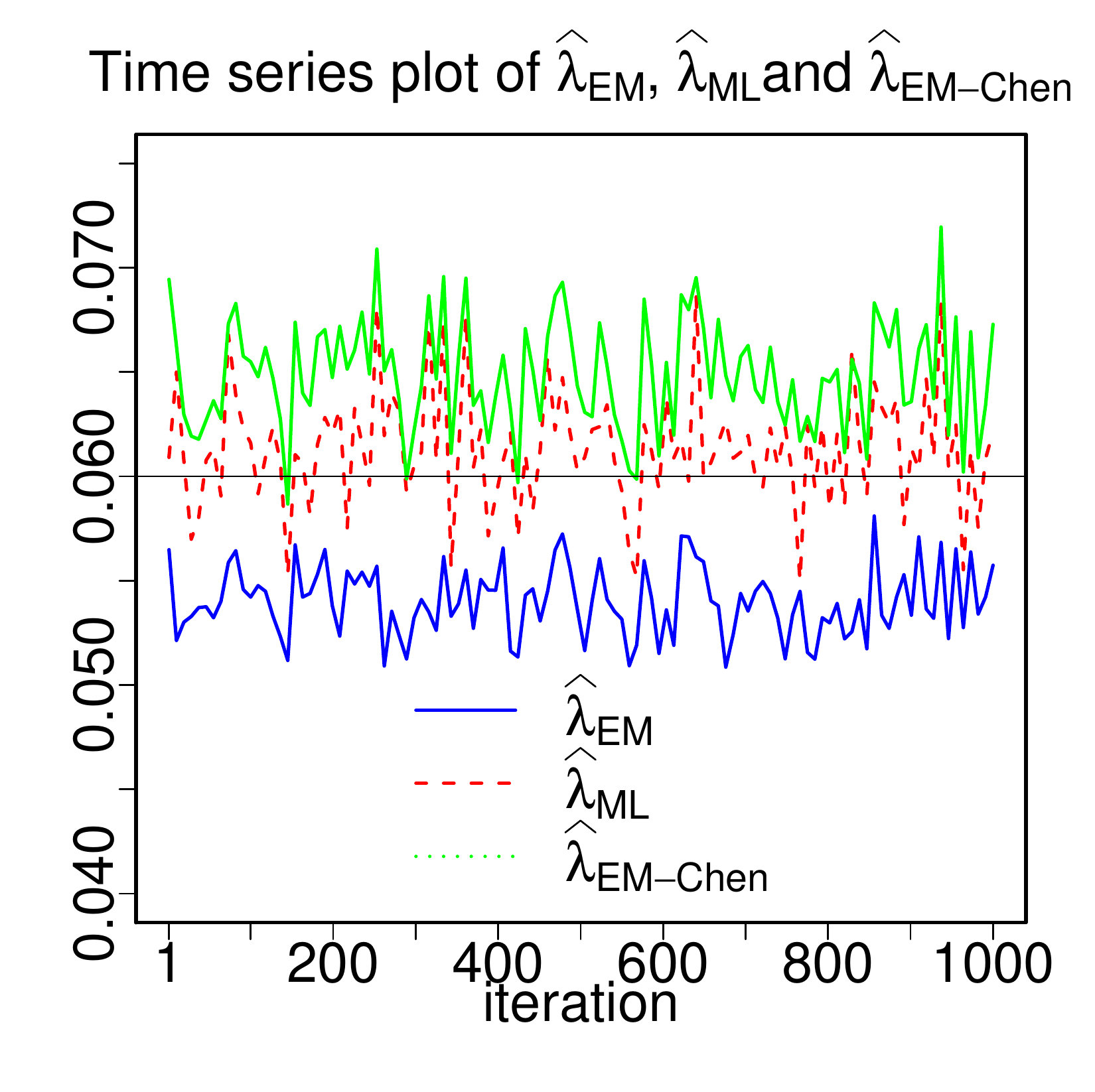}\\
\includegraphics[width=50mm,height=50mm]{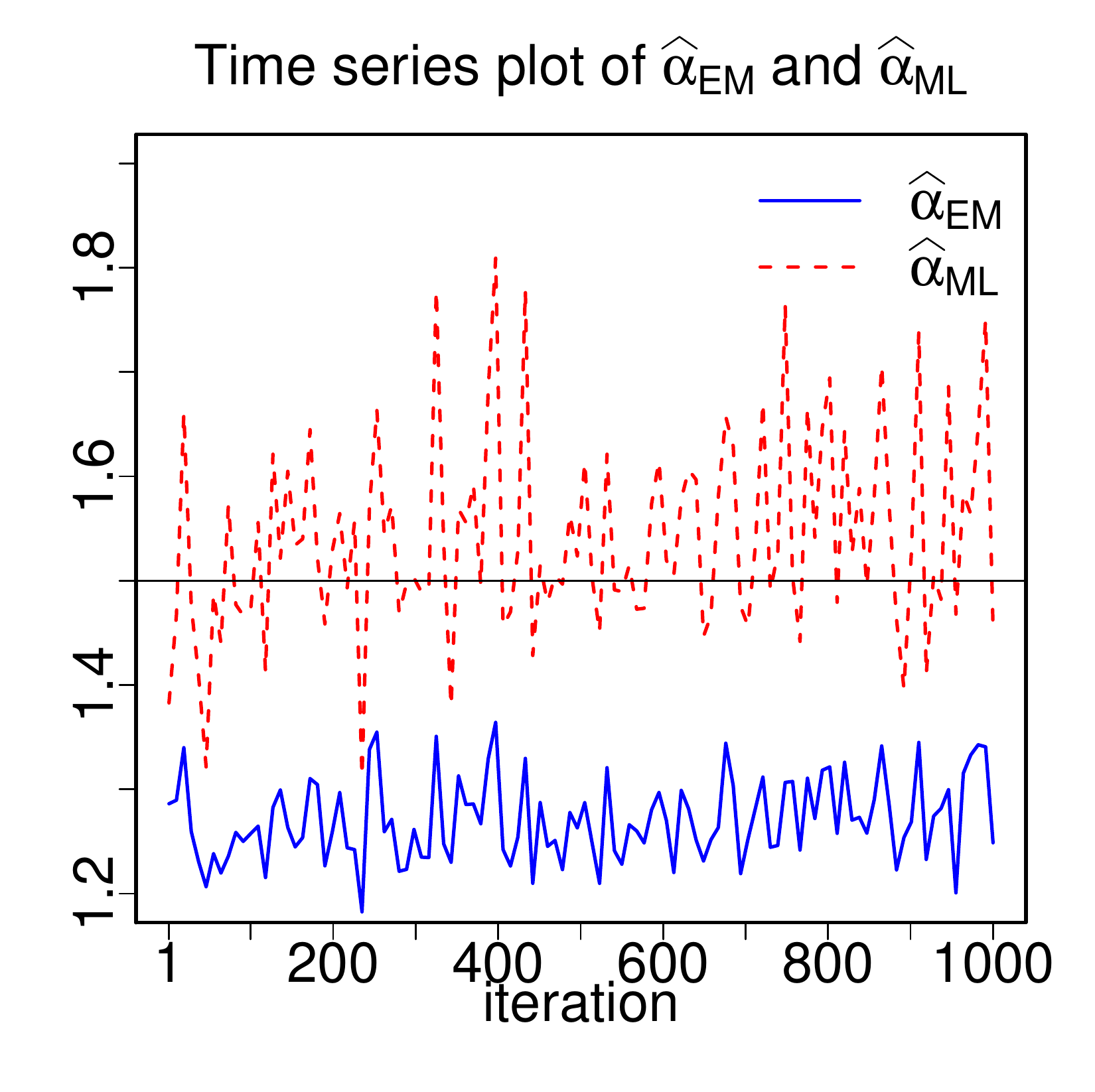}&
\includegraphics[width=50mm,height=50mm]{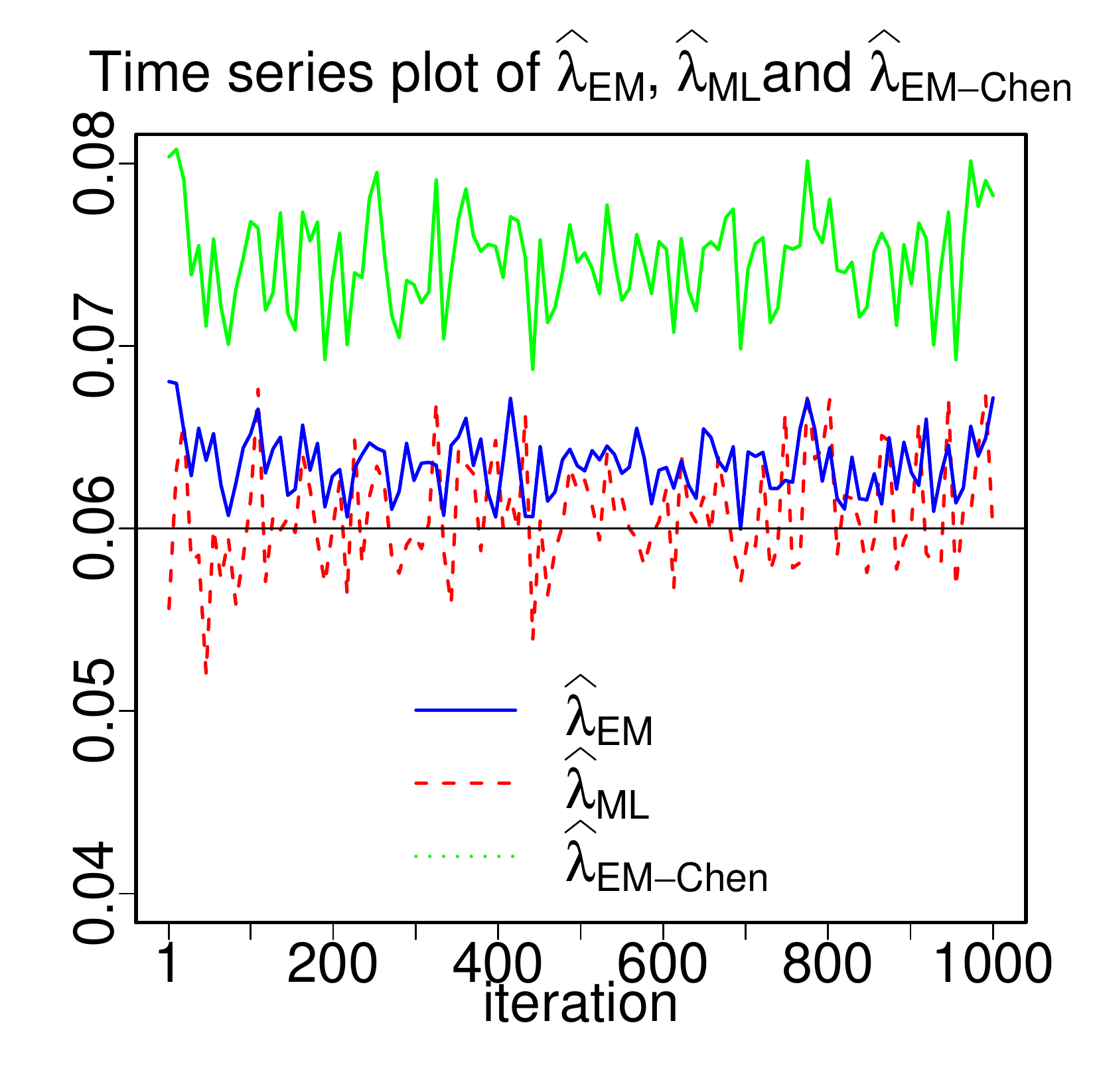}\\
\end{tabular}}
\caption{Time series plot of $\hat{\alpha}_{EM}$, $\hat{\alpha}_{ML}$, and $\hat{\alpha}_{EM-Chen}$ under settings $p_{(3)}$ (top row) and $p_{(4)}$ (bottom row).}
\label{fig2}
\end{figure}

Here, we perform a simulation study to compare the performance of three estimators including: EM algorithm, ML, and EM algorithm proposed by Chen and Lio (2010) for estimating the parameters of GE distribution when items lie under progressive type-I censoring scheme. For simulating a $\{X_i,R_i,T_i\}_{i=1}^{m}$ scheme we use the algorithm proposed by Chen and Lio (2010). We consider four scenarios as:
\begin{align}
p_{(1)}&=(0.25,0.25,0.25,0.25,0.50,0.50,0.50,0.50,1),\nonumber\\
p_{(2)}&=(0.50,0.50,0.50,0.50,0.25,0.25,0.25,0.25,1),\nonumber\\
p_{(3)}&=(0,0,0,0,0,0,0,0,1),\nonumber
\end{align}
and $p_{(4)}=(0.25,0,0,0,0,0,0,0,1)$. Under each of above four scenarios, we simulate $n=112$ observations from GE distribution with shape parameter $\alpha=1.5$ and rate parameter $\lambda=0.06$ and $m=9$ pre-specified inspection times including: $t_1=5.5$, $t_2=10.5$, $t_3=15.5$, $t_4=20.5$, $t_5=25.5$, $t_6=30.5$, $t_7=40.5$, $t_8=50.5$, and termination time is $t_9=60.5$. These settings was used by Chen and Lio (2010). We run simulations for 1000 times when the ML method, proposed EM algorithm in this paper (called here EM), and proposed EM algorithm by Chen and Lio (2010) (called here EM-Chen) take part in the competition. We note that the starting values for implementing both of EM and EM-Chen algorithms are $\alpha^{(0)}=1$ and $\lambda^{(0)}=0.5$. The stopping criterion for both algorithms is $\max\bigl\{\big|\alpha^{(t+1)}-\alpha^{(t)}\big|,\big|\lambda^{(t+1)}-\lambda^{(t)}\big|\bigr\}\leq0.000001$; for $t=0,\dots,100$. The time series plots of the estimators are displayed in Figures (\ref{fig1})-(\ref{fig2}). The summary statistics including bias and mean of squared errors (MSE) of estimators are given in Table \ref{tab1}. Recall that the EM and EM-Chen algorithms give the same estimators for the shape parameter and hence time series plot of $\hat{\alpha}_{EM-Chen}$ disappeared in left-hand side subfigures of Figures (\ref{fig1})-(\ref{fig2}). As it is seen from Table \ref{tab1}, proposed EM algorithm outperforms EM-Chen algorithm under the first, second, and fourth scenarios in terms of bias, and it outperforms the EM-Chen algorithm in all four scenarios in the sense of MSE. Also, the EM algorithm shows better performance than the ML approach under the first scenario in the sense of both bias and MSE criteria.
\begin{table}[h]
\center
\caption{Bias and MSE of $\hat{\alpha}_{EM}$, $\hat{\lambda}_{EM}$, $\hat{\alpha}_{ML}$, $\hat{\lambda}_{ML}$, $\hat{\alpha}_{EM-Chen}$, and $\hat{\lambda}_{EM-Chen}$ under four settings $p_{(1)}$, $p_{(2)}$, $p_{(3)}$, and $p_{(4)}$.} 
\begin{tabular}{cccccc} 
\cline{1-6} 
scenario&Estimator&bias $\hat{\alpha}$& MSE $\hat{\alpha}$ & bias $\hat{\lambda}$ & MSE $\hat{\lambda}$\\ \cline{1-6}
\multirow{2}{*}{$p_{(1)}$}& EM &-0.03470&0.03709&0.01747&0.00033\\
&ML &0.05680&0.10186&0.00119&0.00012\\ 
&EM-Chen &-0.03470&0.03709&0.04212&0.00203\\ \cline{1-6}
\multirow{2}{*}{$p_{(2)}$}& EM &0.10301&0.05822&0.03990&0.00162\\
&ML &0.07546& 0.16765& 0.00222& 0.00027\\ 
&EM-Chen &0.10301&0.05822&0.08084&0.00701\\ \cline{1-6}
\multirow{2}{*}{$p_{(3)}$}& EM &-0.21885&0.05793& -0.00581&0.00005\\ 
&ML &0.05504&0.06842&0.00140&0.00006\\ 
&EM-Chen &-0.21885&0.05793&0.00518&0.00007\\ \cline{1-6}
\multirow{2}{*}{$p_{(4)}$}& EM &-0.23154&0.06314&0.00364&0.00003\\ 
&ML & 0.05017&0.06794&0.00101&0.00007\\ 
&EM-Chen &-0.23154&0.06314&0.01484&0.00027\\ \cline{1-6}
\end{tabular} 
\label{tab1}
\end{table} 
\section{Conclusion}
We have discovered that the EM algorithm proposed by Chen and Lio (Computational Statistics and Data Analysis 54: 1581-1591, 2010) for estimating the parameters of generalized exponential distribution under progressive type-I censoring scheme is incorrect. Here, the corrected EM algorithm is proposed and then a comparison study have been made to discover differences. Theoretically there is no difference between shape estimators of our proposed EM algorithm and that proposed by Chen and Lio (2010). However, for the rate parameter the difference is quite significant. A simulation study have been performed to show visually the differences between performance of our proposed EM algorithm, maximum likelihood estimators, and EM algorithm proposed by Chen and Lio (2010). We note that both of our proposed EM algorithm and EM algorithm proposed by Chen and Lio (2010) converge under all four scenarios before 20 iterations.

\end{document}